\documentclass[twocolumn]{article}

\usepackage{graphicx,latexsym,amsmath,amssymb}

\newtheorem{theorem}{Theorem}

\newtheorem{lemma}{Lemma}

\def\dotqi{\dot q^{\hbox{\,}i}}

\def\dotqj{\dot q^{\hbox{\,}j}}

\def\hatV{\widehat V}
\def\hatC{\widehat C}

\def\wg{\widehat g}
\def\wC{\widehat C}
\def\wV{\widehat V}

\begin{document}

\title{Matching, linear systems, and the ball and beam}

\author{F. Andreev${}^{1,3}$ , 
D. Auckly${}^{1,4}$,               
L. Kapitanski ${}^{1, 2,4}$,   
 S. Gosavi${}^{1,5}$, 
W. White${}^{1,5}$, 
 A. Kelkar${}^{1,6}$ } 
\date{${}^3$ Department of Mathematics, 
        Western Illinois University, 
         Macomb, IL 61455, USA \\
${}^4$ Department of Mathematics, 
                Kansas State University,  
                Manhattan, KS 66506, USA \\
${}^5$ Department of Mechanical and Nuclear Engineering,  
        Kansas State University,  Manhattan, KS 66506, USA \\
${}^6$ Department of Mechanical Engineering, 
Iowa State University, Ames, IA 50011, USA}
\maketitle
\setcounter{footnote}{1}
\footnotetext{Supported in part by NSF grant CMS 9813182 }
\setcounter{footnote}{2}
\footnotetext{Supported in part by NSF grant DMS 9970638}

\begin{abstract}
A recent approach to the control of underactuated systems is to
look for control laws which will induce some specified structure
on the closed loop system.  In this 
paper, we describe one matching condition and
an approach for finding all control laws that fit the condition. After
an analysis of the resulting control laws for linear systems, we 
present the results from an experiment on a nonlinear ball and beam system. 
\end{abstract}

\section{ Underactuated systems and the matching condition}

Over the past five years several researchers have proposed nonlinear 
control laws for which the closed loop system assumes some special form, 
see the controlled Lagrangian method of  \cite{BLM1,BLM3,BLM5}
the generalized matching conditions of \cite{H1,H2,H3},
the 
interconnection and damping assignment passivity based control of \cite{BOS}, 
the $\lambda$-method of \cite{AKW,AK2}, and the references therein. 
In this paper we describe the implementation of the $\lambda$-method of \cite{AKW} on  
a ball and beam system. For the readers convenience we start with the statement of 
the main theorem on $\lambda$-method matching control laws (Theorem 1). We 
also present an indicial derivation of the main equations.  We then prove a new theorem 
showing that the family of matching control laws of any linear time invariant  
system contains all linear state feedback control laws (Theorem 2). 
We next present the general solution 
of the matching equations for the Quanser ball and beam system. (Note, that this 
system is different from the system analyzed by Hamberg, \cite{H1}.)  
As always, the general solution contains several free functional parameters 
that may be used as tuning parameters. We chose these arbitrary functions 
in order to have a fair comparison with the manifacturer's linear control law.
Our laboratory tests confirm the predicted stabilization. 
This was our first experimental test of the $\lambda$-method. We later tested 
this method on an inverted pendulum cart, \cite{AAKKW2}.

Consider a system of the form
\begin{equation}
 g_{rj}\ddot x^j\,+\, [j\,k,\,r]\,\dot x^j\,\dot x^k\,+\,
 C_{r}\,+\,\frac{\partial V}{\partial x^r}\,=\,u_r\,,
\label{eq}
\end{equation}
\noindent $r=1,\dots, n$, where $g_{ij}$ denotes the mass-matrix, 
$C_r$ the dissipation, $V$ the potential energy, 
$[i\,j,\,k]$  the Christoffel symbol of the first kind, 
\begin{equation}\label{christ1}
  [jk,i]= \frac{1}{2}
  \left(\frac{\partial  g_{ij}}{\partial q^k} +
  \frac{\partial  g_{ki}}{\partial q^j}
    - \frac{\partial  g_{jk}}{\partial q^i} \right), 
\end{equation} 
 and $u_r$ is the applied actuation. 
To encode the fact that some degrees of freedom are unactuated, 
the applied forces and/or torques are restricted to satisfy 
$P^i_jg^{jk}u_k=0$, where $P^i_j$ is a $g$-orthogonal projection. 
The matching conditions come from this restriction together with the requirement 
that the closed loop system takes the form
$$
 \wg_{rj}\ddot x^j
+ \widehat {[j\,k,\,r]}\dot x^j\dot x^k+
 \wC_{r}+\frac{\partial  \wV}{\partial x^r}=0,
\, r=1,\dots, n,
$$
for some choice of $\widehat g$, $\widehat C$, and $\widehat V$. 
The matching conditions read
\begin{equation}\label{match1}
\begin{aligned}
 & P^r_k \left( \!\Gamma_{ij}^k -  \widehat \Gamma_{ij}^k \!\right)=0,\, 
  P^r_k \left(\! g^{ki}C_i - \wg^{ki} \wC_i \!\right) =0, \\   
  & P^r_k \left(\! g^{ki}\frac{\partial V}{\partial q^i} 
   - \widehat g^{ki}\frac{\partial \hatV}{\partial q^i}\!\right)
     =0,  
\end{aligned}
\end{equation}
where $\Gamma_{ij}^k$ is the Christoffel symbol of the second kind,
\begin{equation}\label{christ2}
\Gamma_{ij}^k=g^{k\ell}[ij,\ell]. 
\end{equation} 
If the matching conditions (\ref{match1}) hold, the control law will be given by 
\begin{equation}\label{control1}
\begin{aligned}
  u_r & = g_{rk}  (\Gamma_{ij}^k - \widehat \Gamma_{ij}^k) \dotqi \dotqj  +
\left(C_r-\wC_r \right) \\
 & +g_{rk}\left( g^{ki}\frac{\partial V}{\partial q^i}    - \widehat g^{ki}\frac{\partial \hatV}{\partial q^i}\right).
\end{aligned}
\end{equation}
The motivation for this method is that 
$\widehat H = \frac{1}{2}\widehat g_{ij} \dotqi\dotqj + \widehat V$
is a natural candidate for a Lyapunov function because 
$\frac{d}{dt} \widehat H =-\widehat g_{ij}\widehat c^i  \dotqj$.
Following \cite{AKW}, introduce new
variables $ \lambda^k_i=g_{ij}\widehat g^{jk}$. We have
\begin{theorem}
The functions $\widehat g_{ij}$, $\wV$, and $\wC$ satisfy (3) 
in a neighborhood of $x_0$ if and only if 
\begin{equation}
\begin{aligned}
P^r_k \left( g^{ki}C_i - \wg^{ki} \wC_i \right) &=0, \\  
P^r_k \left( g^{ki}\frac{\partial V}{\partial q^i} 
   - \widehat g^{ki}\frac{\partial \hatV}{\partial q^i}\right)
     &=0,
\end{aligned}\notag
\end{equation} 
and the following conditions hold. 
First, there exists a hypersurface containing $x_0$ and transverse to 
each of the vectorfields $\lambda^\ell_i P^i_j \partial/\partial x^\ell$ 
on which $\wg_{ij}$ is invertible and symmetric  and satisfies  
$$ 
g_{ki}P^k_\ell=\lambda^j_k P^k_\ell\wg_{ji}. 
$$
Second,  
$\lambda^i_jP^j_k$ and $\widehat g_{ij}$ satisfy 
\begin{gather}
P^s_k P^r_t \left(g_{\ell s}\frac{\partial \lambda^\ell_r}{\partial q^j}
  +[\ell j,s]\lambda^\ell_r
 -[rj,i]\lambda^i_s \right.  \notag\\ 
\left. +g_{ir}\frac{\partial \lambda^i_s}{\partial q^j}
  +[ij,r]\lambda^i_s
 -[sj,\ell]\lambda^\ell_r \right)
     =0, \label{one}\\ 
\lambda^\ell_rP^r_t\frac{\partial \widehat g_{nm}}{\partial q^\ell}
  +\widehat g_{\ell n}\frac{\partial (\lambda^\ell_rP^r_t)}{\partial q^m}
  +\widehat g_{\ell m}\frac{\partial (\lambda^\ell_rP^r_t)}{\partial q^n}  \notag\\ 
= 
P^\ell_t\frac{partial g_{nm}}{\partial q^\ell}+
  + g_{\ell n}\frac{\partial P^\ell_t}{\partial q^m}
g_{\ell m}\frac{\partial P^\ell_t}{\partial q^n}.
\end{gather}
\end{theorem}
Although the proof of this proposition may be found in 
\cite{AKW}, \cite{AK1}, and \cite{AK2},  
for convenience, we include an indicial notation 
derivation of equations (\ref{one}) and (7). 
Substitute equations (\ref{christ1}),  (\ref{christ2}) for both $\Gamma^k_{ij}$ and 
$\widehat \Gamma^k_{ij}$ 
into the first of equations (\ref{match1}) 
 and multiply the result by the scalar $2$ to obtain:
$$
\begin{aligned}
  P^r_k &
     \widehat g^{k \ell}
     \frac{\partial \widehat g_{ij}}{\partial q^\ell}
- P^r_k \widehat g^{k \ell}\frac{\partial \widehat g_{\ell i}}{\partial q^j}
-P^r_k \widehat g^{k \ell} \frac{\partial \widehat g_{j\ell}}{\partial q^i} \\
 & = P^r_k\,g^{k \ell}
   \frac{\partial g_{ij}}{\partial q^\ell}  
       -P^r_k\,g^{k \ell}\frac{\partial  g_{\ell i}}{\partial q^j}
     -P^r_k\,g^{k \ell} \frac{\partial g_{j\ell}}{\partial q^i} \,.
\end{aligned}
$$
Multiply by $g_{rt}$ and use that $P$ is self-adjoint i.e.,
$P^k_ig_{kj}=g_{ik}P^k_j$, to get
\begin{equation}\label{eight}
\begin{aligned}
  P^r_t &
    \lambda^\ell_r \frac{\partial \widehat g_{ij}}{\partial q^\ell} 
 - P^r_t   \lambda^\ell_r  \frac{\partial \widehat g_{\ell i}}{\partial q^j}
 - P^r_t   \lambda^\ell_r  \frac{\partial \widehat g_{j\ell}}{\partial q^i}
 \\ 
  & = 
P^r_t \frac{\partial g_{ij}}{\partial q^r} 
- P^r_t \frac{\partial  g_{r i}}{\partial q^j}
-P^r_t \frac{\partial g_{jr}}{\partial q^i}\;.
\end{aligned}
\end{equation}
Use  $\;P^r_t   \lambda^\ell_r  \frac{\partial \widehat g_{\ell i}}{\partial q^j}
= \frac{\partial (P^r_t   \lambda^\ell_r\widehat g_{\ell i})}{\partial q^j}-
\widehat g_{\ell i}\frac{(\partial P^r_t   \lambda^\ell_r)}{\partial q^j}\;$  and 
\begin{equation}\label{lam}
\lambda^\ell_r\wg_{\ell i}=g_{ri}
\end{equation}
in (\ref{eight}) to obtain (7). 
To derive (\ref{one}), first, differentiate 
(\ref{lam}) 
with respect to 
$q^j$ to get 
\begin{equation}\label{thirt}
   \lambda^\ell_r \frac{\partial \widehat g_{\ell i}}{\partial q^j}
   = \frac{\partial g_{r i}}{\partial q^j} - \widehat g_{\ell i} 
     \frac{\partial   \lambda^\ell_r}{\partial q^j}.
\end{equation}
Substitute equation (\ref{thirt}) into equation (\ref{eight}) and obtain
\begin{equation} \label{shit}
\begin{aligned}
 & P^r_t \left(
\widehat g_{\ell i}\frac{\partial \lambda^\ell_r}{\partial q^j}
+\widehat g_{\ell j}\frac{\partial \lambda^\ell_r}{\partial q^i}
+\lambda^\ell_r\frac{\partial \widehat g_{ij}}{\partial q^\ell}\right) \\ 
& = P^r_t  \left(
    \frac{\partial g_{ri}}{\partial q^j}
      +\frac{\partial g_{rj}}{\partial q^i}
   - \frac{\partial  g_{r i}}{\partial q^j}
     -\frac{\partial g_{jr}}{\partial q^i}
     +\frac{\partial g_{ij}}{\partial q^r} \right)\,.
\end{aligned}
\end{equation}
Multiply by $-P^s_k \lambda^i_s$,   use (\ref{lam}) and (\ref{thirt}) to obtain
\begin{equation}\label{twelve}
\begin{aligned}
 & P^s_k P^r_t \left(
  g_{\ell s}\frac{\partial \lambda^\ell_r}{\partial q^j}
 +\lambda^\ell_r\frac{\partial g_{js}}{\partial q^\ell}
 -\lambda^i_s\frac{\partial g_{ij}}{\partial q^r} \right) \\ 
 & = P^s_k P^r_t \left( 
\widehat g_{ij}
\lambda^\ell_r\frac{\partial \lambda^i_s}{\partial q^\ell}
- \lambda^i_s\widehat g_{\ell j}\frac{\partial \lambda^\ell_r}{\partial q^i}\right) .
\end{aligned}
\end{equation}
Finally, to obtain (\ref{one}), 
add to equation (\ref{twelve}) 
an equation  obtained from (\ref{twelve}) 
by interchanging $k$ and $t$, $r$ and $s$, $\ell$
and $i$.

\section{ Matching and constant coefficient linear systems}

In this section,  we prove that for  linear
time invariant systems any linear 
full state feedback control law is a solution to the matching equations. 
\begin{theorem}  
 When applied to linear, time-independent  
systems, the family of matching control laws contains all linear state feedback laws.
\end{theorem}

Choose coordinates $q^i$ so that the desired equilibrium is at the origin,  
$V=V_{ij}q^iq^j+v_k q^k$, and $C_i=C_{ij}\dotqj$, 
where $g_{ij}$, $V_{ij}$, $v_k$,
$C_{ij}$, and $P^r_k$ are constant, and $P^r_k$  
has rank $n_u$. 
Clearly,  $\widehat g_{ij}$, $\widehat V=\widehat V_{ij}q^iq^j$, 
and $\widehat C_i=\widehat C_{ij}\dotqj$ is a solution to the matching 
equations when $\widehat g_{ij}$, 
$\widehat V_{ij}$, and
$\widehat C_{ij}$ are constant provided $\widehat g_{ij}$ and 
$\widehat V_{ij}$ are symmetric,  
$P^r_k \left( g^{ki} V_{ij} - \widehat g^{ki} \widehat V_{ij}\right) =0$,
and
$P^r_k \left( g^{ki} C_{ij} - \widehat g^{ki} \hatC_{ij}\right) =0$.  
Let $u_k = v_k+a_{ki}q^i+b_{ki}\dotqi$ 
be an arbitrary
linear control law, satisfying $P^r_kg^{k\ell}u_\ell=0$. 
Comparison  with equation (3) gives
$$
g_{rk}\left( g^{ki} V_{ij} - \widehat g^{ki} \hatV_{ij}\right) = a_{rj}, 
$$  
and
$$
g_{rk}\left( g^{ki} C_{ij} - \widehat g^{ki} \hatC_{ij}\right) =b_{rj}.
$$ 
Thus, 
$$
\widehat V_{\ell j} = \widehat g_{\ell p}g^{pr}(V_{rj}-a_{rj})
$$ 
and 
$$ 
\widehat C_{\ell j} = \widehat g_{\ell p}g^{pr}(C_{rj}-b_{rj}).
$$  
It remains to check that we can find a symmetric, nondegenerate 
matrix $\widehat g^{ki}$ so that the resulting $\widehat V_{\ell j}$ is also symmetric. The symmetry of $\widehat V_{\ell j}$ will follow if we have
$$ 
\widehat g_{\ell p}g^{pr}(V_{rj}-a_{rj})-
 \widehat g_{j p}g^{pr}(V_{r\ell}-a_{r\ell})=0,
$$
and, 
therefore, we need to find a symmetric, nondegenerate  
matrix $\widehat g_{\ell p}$ satisfying 
this equation. The existence of such matrix is guaranteed by the following simple 
observation.
\begin{lemma}
Given any real $n\times n$ matrix $R$, there is a nondegenerate 
symmetric matrix $X$ so that 
$$
RX-X^T R^T=0.
$$
\end{lemma}
Indeed, 
setting $X=QYQ^T$, results in the following equation for $Y$:
$$ 
Q^{-1}RQY-Y^T(Q^{-1}RQ)^T=0.
$$ 
Hence, 
without loss of generality we may assume that $Q^{-1}RQ$ is a real 
Jordan block (see \cite{HJ}), 
$$
\begin{pmatrix}
\lambda & 1 & 0  & \dots \\ 
0 & \lambda & 1  & \dots \\
 & & \dots &   
\end{pmatrix}
$$ 
or 
$$
\begin{pmatrix}
  \begin{matrix} a & -b\\ b & \phantom{-}a \end{matrix}
   &
  \begin{matrix} 1 & \phantom{-}0\\ 0 & \phantom{-}1 \end{matrix} 
   &
   \begin{matrix} 0 & \phantom{-}0\\ 0 & \phantom{-}0 \end{matrix} 
   & \dots \\
   \begin{matrix} 0 & \phantom{-}0\\ 0 & \phantom{-}0 \end{matrix}
   & \begin{matrix} a & -b\\ b & \phantom{-}a \end{matrix}
   &
  \begin{matrix} 1 & \phantom{-}0\\ 0 & \phantom{-}1 \end{matrix} 
   & \dots \\ 
   \dots & \dots & \dots & \dots 
\end{pmatrix}\;.
$$ 
In each case 
$
Y=\begin{pmatrix}
0 &  \dots & 0 & 1\\ 
0 &   \dots & 1 & 0\\
 & & \dots & & \\ 
1 & 0 & \dots &  0 
\end{pmatrix}
$ 
solves the equation. 

Note that the result of Lemma 1 is true for matrices with coefficients in 
any field. This is proved in \cite{TZ}.

\section{Example: The Ball and Beam}

In order to demonstrate the approach described above, we have implemented
one of the control laws from the family of control laws described in the first section 
on  a ball and beam system, Figure 1 (this
 system is commercially available from Quanser
 Consulting, Ontario, Canada).

\begin{figure}
\centerline{\includegraphics[width=7.5cm]{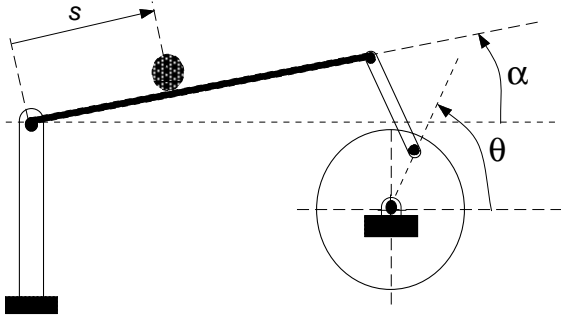}}
\caption{The ball and beam system}\label{fig:bb}
\end{figure}

The $s-$coordinate is unactuated, 
the $\theta-$coordinate is actuated by the servo, and the objective is to bring the ball
to the center of the beam. 
The physical parameters of the system are given in Table 1.
\begin{table*}[htb]
\caption{The physical parameters of the system.}
\label{table:1}
\newcommand{\m}{\hphantom{$-$}}
\newcommand{\cc}[1]{\multicolumn{1}{c}{#1}}
\renewcommand{\tabcolsep}{1.65pc} 
\renewcommand{\arraystretch}{1.2} 
\begin{tabular}{@{}llll}
\hline
 $\ell_b$ = length of the beam &  $=0.43\text{m}$ 
& $I_B =\frac{2}{5} m_B r_B^2 $ = ball 			
inertia &  $=4.25\times 10^{-6}\text{Kg m}^2$  \\
 $\ell_l$  = length of the link  & $=0.11\text{m}$ & $I_b$   = inertia of the beam    &  	
		$=.001\text{Kg m}^2$  \\
 $r_g$     = radius of the gear  &    $=0.03\text{m}$ & $I_s$    = effective servo inertia&  
$=0.002\text{Kg m}^2$ \\
 $r_B$     = radius of the ball  &    $=0.01\text{m}$ & $g$    = gravitational acceleration     
&  $=9.8\text{m/}\text{s}^2$\\ 
 $m_B$     = mass of the ball    &    $=0.07Kg$ & $s_0$  = desired equilibrium position  
&  $=0.22\text{m}$  \\
 $m_b$     = mass of the beam    &    $=0.15Kg$  & $c_0$  = inherent servo dissipation
  &  $=9.33\times 10^{-10} \text{Kg m}^2/\text{s}$ \\
 $m_l$     = mass of the link    &    $=0.01Kg$ & & \\
\hline
\end{tabular}
\end{table*}

\noindent One can express $\alpha$ as a concrete function of $\theta$  
from the kinematic relation
$$
\begin{aligned}
 & \left( \ell_b (1-\cos(\alpha)) -r_g (1-\cos(\theta))\right)^2 \\ 
& +
\left( \ell_b \sin(\alpha)+\ell_l-r_g \sin(\theta)\right)^2=\ell_l^2.
\end{aligned}
$$ 
The kinetic energy of the system is 
$$
T=\frac12 m_b s^2\dot\alpha^2+\frac12 I_B(\dot\alpha+\frac{1}{r_B}\dot s)^2
+\frac12 I_b\dot\alpha^2+\frac12 I_s\dot\theta^2.
$$
The potential energy is 
$$
\begin{aligned}
& V=\frac12 m_{l} \;g \;r_g\sin(\theta)
+\frac12 (m_b+m_l)\;g \;\ell_b\sin(\alpha)\\ 
&+
m_B\; g \;s \sin(\alpha)\,,
\end{aligned}
$$
and the dissipation is $C_1=0$, $C_2=c_0\dot\theta$.
After rescaling, we get 
$$
\begin{aligned}
 (1-\cos(\alpha) &-  a_2 (1-\cos(\theta)))^2 \\ 
  &+(\sin(\alpha)+a_1-a_2\sin(\theta))^2=a_1^2\,, 
\end{aligned}
$$
$$
T=  \frac{1}{2}\,\dot s^{2} +
\frac{1}{2}  (a_4 +(a_3+5/2\,s^{2})
(\alpha^\prime(\theta))^2 )\dot \theta^2 
+ \alpha^\prime(\theta)\,\dot s\,\dot \theta \;,
C_1=0\,,\qquad  C_2=a_7\,\dot\theta\,, 
$$
$$
V=a_5\sin(\theta)+(s+a_6)\sin(\alpha(\theta)),
$$
$C_1=0$, and $C_2=a_7\dot\theta$,
where the $a_k$ are 
the dimensionless parameters, 
\begin{gather} 
a_1=\frac{\ell_l}{\ell_b},\quad    a_2=\frac{r_g}{\ell_b},\quad 
a_3=\frac{(I_b+I_B)}{I_B},\quad   a_4=\frac{I_s}{I_B},\notag\\  
a_5=\frac{m_lr_g}{2m_Br_B},\quad  
a_6=\frac{\ell_b(m_b+m_l)}{2m_Br_B},\notag\\  
a_7=\left(\frac {5}{2 \;r_B^3\; g} \right)^\frac12 \,\frac{c_0}{ m_B}\;.\notag  
\end{gather}
The notation $\,{}^\prime\,$ is used to 
denote a derivative of a function of one variable.  
For general underactuated systems, 
the use of the powerful $\lambda$-method to solve the matching equations 
is discussed in \cite{AKW,AK2}.  For systems with two degrees of freedom, 
the $\lambda$-method produces the  general solution 
to the matching equations in an explicit form, \cite{AK1}.  
When applied to the ball and beam system, the explicit family of control laws 
is given by equation (\ref{control1}) with the following expressions for 
$\wg$, $\wV$, and $\wC$, where 
\begin{gather} 
\widehat g_{11}= \psi^2(\alpha)
\left( h(y(s,\theta))+10\,\int_0^\alpha
      \frac{ d\varphi}{\mu_1^{\prime}(\varphi)\psi^2(\varphi)}
\right),\notag\\ 
\widehat g_{12} =\frac{1}{\mu} (g_{11}-\sigma \widehat g_{11}),
\; 
\widehat g_{22} =\frac{1}{\mu} (g_{12}-\sigma \widehat g_{12}),\notag 
\end{gather}
\begin{gather}
\begin{aligned}
 \widehat V(s,\theta)= &w(y)+
5(y+s_0)\int_0^\alpha
\frac{\sin(\varphi)}{\mu_1^\prime(\varphi)\psi(\varphi)}\,d\varphi
 \\ 
&\ \\
-&5\,\int_0^\alpha 
\frac{\sin(\varphi)}{\mu_1^\prime(\varphi)\psi(\varphi)}\,
\int_0^\varphi \psi(\tau)\,d\tau\;d\varphi ,
\end{aligned} \notag\\ 
\ \notag\\ 
\wC_1 = (g_{1i}\wg^{i1})^{-1}\left(C_1-g_{1j}\wg^{j2} \wC_2\right), \notag
\end{gather}
where 
\begin{gather}
\mu(s,\theta)=\frac{\mu_1^{\prime}(\alpha(\theta))}{ 5 s\,g_{12}},\notag\\ 
\sigma(s,\theta)=\mu_1(\alpha)-\frac{1}{ 5s}\,\mu_1^{\prime}(\alpha),\notag  \\ 
y\,=\,\psi(\alpha)\,s\,-\,s_0\,+\,\int_0^\alpha \psi(\tau)\,d\tau,
\frac{\mu_1(\kappa)}{\mu_1^{\prime}(\kappa)}\,d\kappa\}.\notag
\end{gather}
\noindent 
Here $h(y)$, $w(y)$, $\mu_1(\alpha)$ are arbitrary functions of one variable, 
and $\wC_2$ is an arbitrary function which is odd in velocities. 

\section{Experimental Results}

Our experiments were conducted on the  Quanser ball and beam system.  
The control signal is a voltage supplied to the servo 
and the sensed output of the system is $s$ and $\theta$ sampled at $300$ Hz.
A Quanser 
MULTIQ$^{\circledR}$ data acquisition card is used for the analog signal 
input and output. 
The velocities are computed via numerical differentiation using the forward 
difference algorithm. The control law produces a voltage signal 
and is supplied through the D/A converter  
to the DC servomotor via an amplifier. The relation between the control voltage,
$v_{in}$, and the torque, $u$ ($=u_2$ in equation (\ref{control1})), is 
$ K_m^2 N_g^2 \dot\theta$, where 
$R_m=$ armature resistance    $=2.6 \;\Omega$, 
 $N_g=$ gear ratio    $=70.5$, 
 $K_m=$ motor torque constant  $=0.00767 $ Volt$\cdot$sec.

Any stabilizing linear control law for this system is specified by four constants.
The nonlinear control laws in our family are specified by the four 
arbitrary {\it functions\/}: $\mu_1(\alpha)$, $h(y)$, $w(y)$, and 
$\widehat C_2(s, \theta, \dot s,\dot\theta)$.
We chose 
\begin{equation}
\begin{aligned}
\mu_1(\alpha) & =1.0849\exp(4.7845 \sin(\alpha)) \\
h(y) &=1.1031,\;  w(y)=0.0023 y^2, \\  
\widehat C_2(s, \theta, \dot s,\dot\theta)
& =-\widehat g_{12}\cdot (1+\dot s^2+10\dot \theta^2)(-\mu\dot s+\sigma \dot \theta).
\end{aligned} \notag
\end{equation}
These functions produce the control law, $u$, in rescaled units. 
The values of the constants $a_1$ through $a_7$ are as follows
\begin{center}
\begin{tabular}{llll}
 $a_1 =\  0.2547$ & $a_5 =\  0.1889$\\
 $a_2 =\  0.0588$ & 
$a_6 =\  42$ \\
 $a_3 =\  236.294$ & 
$a_7=\   5\times 10^{-6}$\\
 $a_4=\  471.126$  & &\\
\end{tabular}
\end{center}
The final 
control signal is obtained by converting back into MKS units 
and using the formula in the preceding paragraph to get the input voltage.
These choices were made from the following considerations.  
The form of the function $\mu_1$ was chosen to simplify 
the integrals in the expressions for $y$, $\psi$, and $\widehat g_{11}$. 
The form of $\widehat C_2$ was chosen to ensure that 
$\widehat C_1 \dot s + \widehat C_2 \dot\theta $ would be positive 
(for $\widehat H$ to be a Lyapunov function). 
Finally, the coefficients in these functions were chosen so
 that the linearization 
of the nonlinear control law would agree with the linear control law 
provided by the manufacturer. 

Extensive numerical simulations done using  Matlab$^{\circledR}$   
 confirm that the nonlinear control law stabilizes the system. 
The linear control law appears to stabilize 
the system for a wider range of initial conditions than the 
nonlinear control law. This is an empirical observation, not 
a mathematical fact. Finding an adequate mathematical framework 
to compare different control laws is a very interesting unresolved problem, 
see \cite{AK1}. Usually, given two locally stabilizing control laws, 
there exist initial conditions stabilized by one but not by the other. 
For example, one set of physically unrealistic 
initial conditions with a large angular 
velocity $\dot \theta =3.6$ (or $158\;\text{rad}/\text{sec}$ in physical units) 
is stabilized by our nonlinear control law but not by the linear one.

\begin{figure}
\centerline{\includegraphics[width=7.5cm]{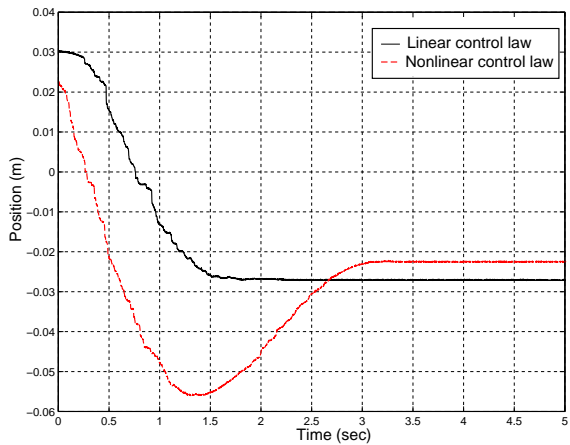}}
\caption{Ball position response}\label{fig:bbxpos}
\end{figure}

\begin{figure}
\centerline{\includegraphics[width=7.5cm]{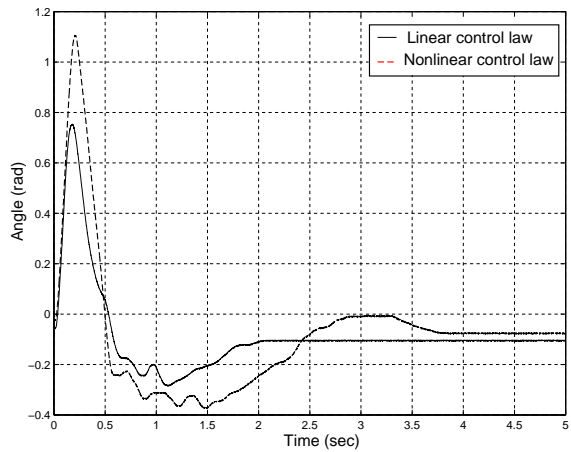}}
\caption{Angular displacement response}\label{fig:bbthpos}
\end{figure}

We have implemented the nonlinear control law in the laboratory. 
The laboratory tests confirm the predicted behavior of the nonlinear 
controller. Figures \ref{fig:bbxpos} and \ref{fig:bbthpos} show a comparison 
of the time histories of the ball position ($s$) and angular displacement ($\theta$) 
for the linear and nonlinear control laws. In both cases 
the control signal reached the saturation limit for a short duration during the 
initial 
rise of the response. The difference in the steady-state values of the 
responses is attributed to a lack of sensitivity of the resistive strip 
used to measure the ball position.  
\section{Conclusions}

The $\lambda$-method produces explicit infinite-dimensional families of 
control laws and simultaneously provides a natural candidate for a 
Lyapunov function. When this method is applied to linear time-invariant  
systems, the resulting family contains all linear state feedback control 
laws (Proposition 2). 
In this paper we also present the results of the first implementation of
a $\lambda$-method matching control law on a concrete physical device, 
the ball and beam system.  
The experimental results agree with theoretical predictions and 
numerical simulations. In our experiments we observe that the linear control law 
performs better than our nonlinear control law for the ball and beam system. 
However, in a later experiment with an inverted pendulum cart, \cite{AAKKW2}, 
we found that a properly tuned  $\lambda$-method matching control law 
performed better than the corresponding linear one.  
At the moment, it is not known for which systems matching control laws 
will perform better. This is an important problem that must be resolved.

\end{document}